# Title: Axiomatic Digital Topology

Author: Vladimir Kovalevsky

Affiliation: Department of Computer Science
University of Applied Sciences Berlin
Luxemburger Str. 10, 13353 Berlin
Germany

Contact: Tel: +49-30-659-60-14; Fax: +49-30-659-08-707
E-mail: kovalev@tfh-berlin.de

Web site: www.kovalevsky.de

# Axiomatic Digital Topology

Vladimir Kovalevsky

**Abstract**

The paper presents a new set of axioms of digital topology, which are easily understandable for application developers. They define a class of locally finite (LF) topological spaces. An important property of LF spaces satisfying the axioms is that the neighborhood relation is antisymmetric and transitive. Therefore any connected and non-trivial LF space is isomorphic to an abstract cell complex. The paper demonstrates that in an *n*-dimensional digital space only those of the (*a*, *b*)-adjacencies commonly used in computer imagery have analogs among the LF spaces, in which *a* and *b* are different and one of the adjacencies is the “maximal” one, corresponding to $3^n-1$ neighbors. Even these (*a*, *b*)-adjacencies have important limitations and drawbacks. The most important one is that they are applicable only to binary images. The way of easily using LF spaces in computer imagery on standard orthogonal grids containing only pixels or voxels and no cells of lower dimensions is suggested.

## 1 Introduction

In 2001 Yung Kong [9] issued the challenge to "construct a simplest possible theory that gives an axiomatic definition of well-behaved 3-d digital spaces". In many publications on digital topology [4, 7] including those by T.Y. Kong and A. Rosenfeld [11], A. Rosenfeld being the author of the graph-based approach to digital topology, the authors interpret “well-behaved” as being in accordance with classical topology with regard to connectedness and validity of the Jordan theorem. One possible way to check the accordance is to adapt the classical axioms of the topology to digital spaces and to compare the connectedness relation in a digital image and in a corresponding topological space. We shall discuss the most important and recent publications concerned with this comparison in Section 5 after we have introduced all necessary definitions.

We prefer to demonstrate here another way since for many practically oriented researchers it is not clear, why the notion of an open subset, that is the key notion of the classical topology, should be applied to digital spaces. Practically oriented researchers may put the following questions: Is it absolutely necessary to use the notion of open subsets satisfying the classical axioms to achieve the purposes of applications? Is it perhaps possible to find quite different axioms and to construct a theory based on these new axioms, which would satisfy all practical demands?

We formulate in Section 2 a new set of axioms, which are “natural” from the point of view of our intuition and of practical demands. Then we prove some theorems and demonstrate the consequences for the theory and for applications of digital spaces, especially for the (*a*, *b*)-adjacency relations commonly used in computer imagery.

We demonstrate in Section 3 that a topological space satisfying our Axioms is a particular case of the classical topological space. In Section 4 we deduce from our Axioms the most important properties of the locally finite spaces satisfying Axioms (ALF spaces) and demonstrate that the neighborhood relation must be antisymmetric and transitive. In Section 5 we discuss the previous work (it is impossible to discuss the previous work before all the necessary notions are defined and explained). In Section 6 we investigate, following the above mentioned challenge by Yung Kong, the pairs of adjacency relations commonly used in computer imagery and demonstrate that in spaces of any dimension *n* only those pairs (*a*, *b*) of adjacencies are consistent, in which exactly one of the adjacencies is the “maximal” one, corresponding to $3^n-1$ neighbors. We also introduce the notion of homogeneous completely connected spaces and demonstrate that the only such space of dimension 3 is the complex whose principal cells are isomorphic to truncated octahedrons, i.e. to polyhedrons with 14 faces. Section 7 is devoted to recommendations for applications. We suggest there techniques which allow one to easily apply the concept of ALF spaces to standard orthogonal grids containing only pixels or voxels and no cells of lower dimensions. The list of abbreviations is to be found in Appendix 2.

## 2 Axioms of Digital Topology

Let us start with an attempt to suggest certain axioms and definitions of topological notions, which are different from those usually to be found in text books of topology, while being comprehensible for practically oriented researchers.

A digital space must obviously be a so called *locally finite space* to be explicitly representable in the computer. In such a space each element has a neighborhood consisting of *finitely many* elements. We don't call the elements "points" since, as we will see, a locally finite space *must* possess elements with different topological properties, which must have different notations.

We suggest the following set of axioms concerned with the notions of connectedness and boundary, which are the topological features most important for applications in computer imagery.

**Definition LFS:** (locally finite space) A nonempty set $S$ is called a *locally finite* (LF) space if to each element $e$ of $S$ certain finite subsets of $S$ are assigned as neighborhoods of $e$.

We shall consider in what follows a particular case of an LF space satisfying the following Axioms 1 to 4 and we shall denote it an ALF space.

**Axiom 1:** For each space element $e$ of the space $S$ there are certain subsets containing $e$, which are neighborhoods of $e$. The intersection of two neighborhoods of $e$ is again a neighborhood of $e$.

Since the space is locally finite, there exists the smallest neighborhood of $e$ that is the intersection of all neighborhoods of $e$. Thus each neighborhood of $e$ contains its smallest neighborhood. We shall denote the smallest neighborhood of $e$ by SN($e$).

**Axiom 2:** There are space elements, which have in their SN more than one element.

**Definition IN:** (incidence) If $b \in$SN($a$) or $a \in$SN($b$), then the elements $a$ and $b$ are called *incident* to each other.

According to the above definition, the incidence relation is symmetric. It is reflexive since $a \in$SN($a$). The notion of incident elements seems perhaps to be similar to the adjacency introduced in [20]. There is, however, an important difference between them because we do not suppose that all elements have the same number of incident elements.

**Definition IP:** (incidence path) Let $T$ be a subset of the space $S$. A sequence ($a_1$, $a_2$, ..., $a_k$), $a_i, \in T$, $i$=1, 2, ..., $k$; in which each two subsequent elements are incident to each other, is called an *incidence path in T* from $a_1$ to $a_k$.

**Definition CN:** (connected) Incident elements are *directly connected*. A subset $T$ of the space $S$ is *connected* iff for any two elements of $T$ there exists an incidence path containing these two elements, which completely lies in $T$.

Let us now formulate axioms related to the notion of a boundary. The classical definition of a boundary (exactly speaking, of the topological boundary or of the frontier) is as follows:

**Definition FR:** (frontier) The topological boundary, also called the *frontier*, of a non-empty subset $T$ of the space $S$ is the set of all elements $e$ of $S$, such that *each* neighborhood of $e$ contains elements of both $T$ and its complement $S–T$.

We shall denote the frontier of $T \subseteq S$ by Fr($T$, $S$).

In the case of a locally finite space it is obviously possible to replace "each neighborhood" by "smallest neighborhood" since according to Axiom 1 each neighborhood of *a* contains the smallest neighborhood of *a*. Now let us introduce the notions of a *thick* and a *thin* frontier.

**Definition NR:** (neighborhood relation) The *neighborhood relation N* is a binary relation in the set of the elements of the space *S*. The ordered pair (*a*, *b*) is in *N* iff $a \in \mathrm{SN}(b)$.

We also write *aNb* for (*a*, *b*) in *N*.

**Definition OT:** (opponents) A pair (*a*, *b*) of elements of the frontier Fr(*T*, *S*) of a subset $T \subset S$ are *opponents* of each other, if *a* belongs to SN(*b*), *b* belongs to SN(*a*), one of them belongs to *T* and the other one to its complement *S*–*T*.

**Definition TF:** (thick frontier) The frontier Fr(*T, S)* of a subset *T* of a space *S* is called *thick* if it contains at least one pair of opponents. Otherwise the frontier is called *thin*.

To justify of the notation "thick" let us remark, that at locations where there are opponent pairs in the frontier, the frontier is doubled: there are two subsets of the frontier, which run "parallel" to each other. These subsets are called border and coborder in [7]. Fig. 1.1 shows some examples.

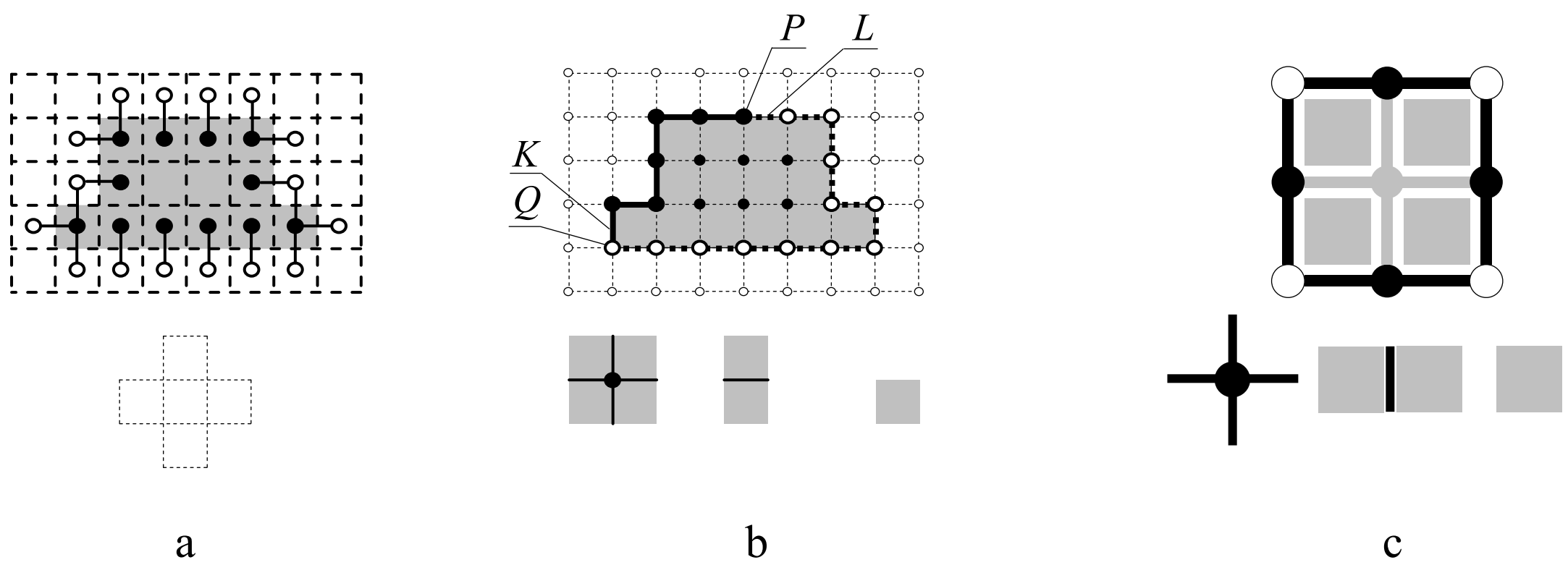

**Fig. 1.1.** Examples of frontiers: a thick frontier (a); a thin frontier (b); a frontier with gaps (c)

The smallest neighborhoods of different space element are shown below the grid. They are different in a, b and c. In Fig. 1.1a the space *S* consists of squares only. The subset *T* is the set of gray squares. The symmetric relation *N* is the well-known 4-adjacency. Elements of Fr(*T*, *S*) are labeled by black or white disks. Squares with black disks belong to *T*, while those with white disks belong to *S*–*T*. In Fig. 1.1a each element of the frontier has at least one opponent. Pairs of opponents are shown by connecting lines. The frontier is thick.

In Fig. 1.1b the space consists of points (represented by disks), lines and squares. The frontier Fr(*T*, *S*) is represented by bold lines and big disks. The gray squares, solid lines and black disks belong to *T*; the white squares, dotted lines and white disks belong to *S*–*T*. The solid line $K \in T$ and the white point $Q \in S-T$ belong to different subsets. Nevertheless, they are no opponents since $Q \notin \mathrm{SN}(K)$. The same is true for the pair (*L*, *P*). Thus the frontier is thin.

**Axiom 3:** The frontier Fr(*T*, *S*) of any subset $T \subset S$ is thin.

According to Definition Fr the frontier of *T* is the same as the frontier of its complement *S*–*T*.

Another important property of the frontier is, non-rigorously speaking, that it must have no gaps, which is not the same, as to say, that it must be connected. More precisely, this means, that the frontier of a frontier *F* is the same as *F*. For example, the frontier in Fig. 1.1c has gaps represented by white disks. Let us explain this. Fig. 1.1c shows a space *S* consisting of points, lines and squares. The relation *N* is for this case non-transitive: the SN(*P*) of a point *P* contains some lines incident to *P* but no squares. The SN of a line contains one or two incident squares, while the SN of a square is the square itself. The subset *T* is represented by gray elements. Its frontier Fr(*T*, *S*) consists of black lines and black points since these elements do not belong to *T*, while their SNs intersect *T*. The white points do not belong to *F*=Fr(*T*, *S*) because their SNs do not intersect *T*. These are the gaps. However, Fr(*F*, *S*) contains the white points because their SNs intersect both *F* and its complement (at the points itself). Thus in this case the frontier Fr(*F*, *S*) is different from *F*= Fr(*T*, *S*).

We shall prove below that the frontier Fr(*F*, *S*) is different from *F*= Fr(*T*, *S*) only if the neighborhood relation is non-transitive, which fact is important to demonstrate that the smallest neighborhoods satisfying our Axioms are open subsets of the space.

**Axiom 4:** The frontier of Fr(*T*, *S*) is the same as Fr(*T*, *S*), i.e. Fr(Fr(*T*, *S*), *S*)=Fr(*T*, *S*).

Let us remain the reader the classical axioms of the topology. The topology of a space *S* is defined if a collection of subsets of *S* is declared to be the collection of *open subsets* satisfying the following axioms:

Axiom C1: The entire set *S* and the empty subset ∅ are open.

Axiom C2: The union of any number of open subsets is open.

Axiom C3: The intersection of a finite number of open subsets is open.

Also an additional so-called separation axiom is often formulated in classical topology:

Axiom C4: The space has the separation property.

There are (at least) three versions of the separation property and thus of Axiom C4 (Fig. 1.2) :

Axiom $T_0$: For any two distinct points *x* and *y* there is an open subset containing exactly one of the points [2].

Axiom $T_1$: For any two distinct points *x* and *y* there is an open subset containing *x* but not *y* and another open subset, containing *y* but not *x*.

Axiom $T_2$: For any two distinct points *x* and *y* there are two non-intersecting open subsets containing exactly one of the points.

A space with the separation property $T_2$ is called *Hausdorff space*. The well-known Euclidean space is a Hausdorff space.

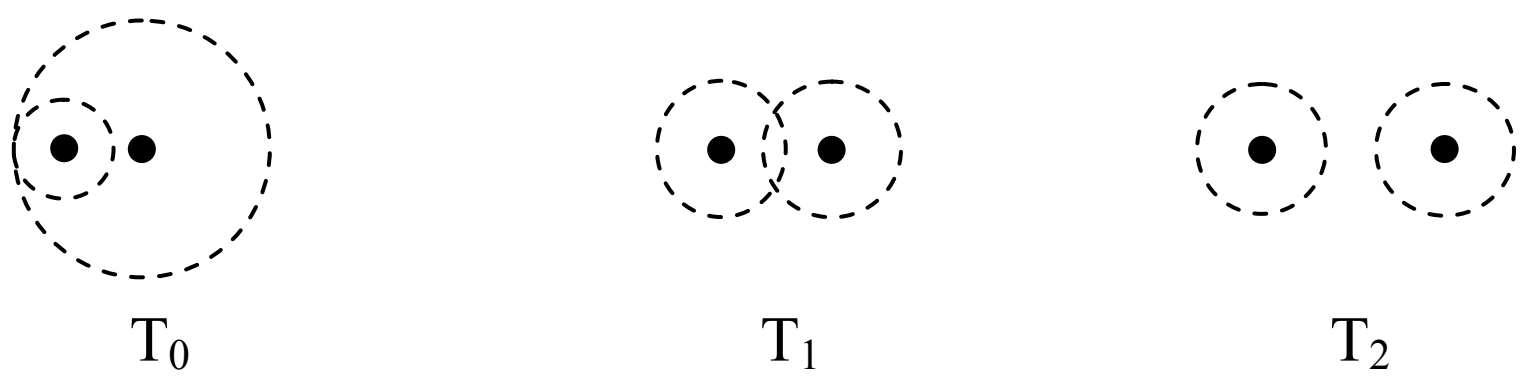

**Fig. 1.2.** A symbolic illustration to the separation axioms

Now we shall demonstrate that the axioms of the classical topology follow as theorems from our set of Axioms 1 to 4.

## 3 Relation between the Suggested and the Classical Axioms

We shall show in the next Section 4 that the classical Axiom $T_0$ follows from the demands that the frontiers must be thin and that the frontier of a subset *T* must be the same as the frontier of its complement $S-T$. For this purpose we shall show now that the neighborhood relation in a space satisfying our Axioms must be antisymmetric.

**Theorem TF:** (thin frontier) An LF space $S$ satisfies the Axiom 3 (Section 2) iff the neighborhood relation $N$ of the space $S$ is *antisymmetric*.

To make the reading easier, we have presented all proofs, except that of the Main Theorem below, in Appendix 1.

There is a possibility to achieve that the frontier be thin for any subset of a space with a *symmetric* neighborhood relation: it is necessary to change the Definition FR of the frontier so that only elements of $T$ may belong to the frontier of $T$ [20]. However, this possibility leads to *different frontiers* of $T$ and of its complement $S-T$, which may be considered as a topological paradox and contradicts the classical definition of the frontier. Thus, the neighborhood relation in a space satisfying our Axioms must be antisymmetric.

Now let us introduce the notion of open subsets of a locally finite space.

**Definition OP:** (open) A subset $O \subset S$ is called *open in S* if it contains no elements of its frontier $\text{Fr}(O, S)$. A subset $C \subset S$ is called *closed in S* if it contains *all* elements of $\text{Fr}(C, S)$.

**Lemma SI:** (SN in subset) A subset $T \subset S$ is open in $S$ according to Definition OP iff it contains together with each element $a \in T$ also its smallest neighborhood $\text{SN}(a)$.

Now let us demonstrate that the first three classical axioms C1 to C3 may be deduced as theorems from our Axioms.

**Theorem OS:** (open subsets) Subsets of an ALF space $S$, which are open in $S$ according to Definition OP, satisfy the classical Axioms $C_1$ to $C_3$ and therefore are open in the classical sense.

In the proof of Theorem OS the classical axioms C1 to C3 of the topology are deduced from our Axiom 1 and Definitions FR and OP. In the following Section we shall demonstrate, that an ALF space is a particular case of a classical topological $T_0$ space.

## 4 Deducing the Properties of ALF Spaces from the Axioms

In the previous Section we have demonstrated that an ALF space is a topological space in the classical sense and that the neighborhood relation is antisymmetric (Theorem TF, Section 3).

Consider the relation "$a \neq b$ and $a \in \text{SN}(b)$". It is usual to call it the *bounding relation*, to denote it by $B$ and to say "$b$ bounds $a$" or "$a$ is bounded by $b$". This notation reflects the fact that $b \in \text{Fr}(\{a\}, S)$. The relation $B$ is irreflexive. According to Theorem TF it must be asymmetric. Now we shall demonstrate that it is transitive.

**Lemma MM:** (minimum and maximum) If $S$ is an ALF space and the bounding relation $B$ is transitive, then $S$ contains elements, which are bounded by no other elements, and it contains elements, which bound no other elements.

The first ones will be called the *minimum* elements and the latter the *maximum* ones.

**Lemma NM:** (no maximum in Fr) Let $T$ be a subset of $S$. If the bounding relation $B$ is transitive, then $\text{Fr}(T, S)$ contains *no maximum elements* of $S$ and for any element $a$ of $S$ the subset $\text{SN}(a)$ contains at least one maximum element of $S$.

**Theorem TR**: (transitive) An LF space satisfies Axiom 4 iff the bounding relation is transitive.

**Corollary HO:** (half-order) The bounding relation *B*, being irreflexive, asymmetric and transitive, is an irreflexive half-order and we can write $a<b$ instead of *aBb*.

**Corollary NO:** (neighborhood is open) The smallest neighborhood of any element *a* of an ALF space is *open* both according to Definition OP and in the classical sense. It is the *smallest open subset* containing *a*.

**Corollary T0:** The smallest neighborhoods in an ALF space satisfy the classical Axiom $T_0$.

The proofs are to be found in Appendix 1.

**Conclusion:** We have demonstrated in Section 3 that the classical axioms C1 to C3 can be deduced as theorems from our Axioms. Now we see that the classical Axiom $T_0$ also follows from our Axioms. This means that an ALF space is a *particular case* of the classical $T_0$ space. We have chosen our Axioms 1 to 4 in the hope that they will be naturally comprehensible for practically oriented researchers. Thus this consideration may serve as an answer to the first question raised in the Introduction: Yes, it is necessary to use the classical open subsets to arrive at practically acceptable notions of connectedness and boundary.

Let us now consider an important particular case of LF spaces known as abstract cell complexes.

**Definition SON:** (smallest open neighborhood) The smallest open subset of the ALF space *S* that contains the element $a \in S$ is called the *smallest open neighborhood* of *a* in *S* and is denoted by SON(*a*, *S*). According to Corollary NO SON(*a*, *S*)=SN(*a*).

**Definition FC:** (face) A space element *a* is called a *face* of the element *b* if $b \in$SON(*a*, *S*). If *a*=*b*, then *a* is a non-proper face of *b*. The *face relation* is reflexive, antisymmetric and transitive. Thus, it is a reflexive partial order in *S* and it is usual to denote it by $a \le b$.

According to Corollary NO the neighborhood relation *N* (Section 2) in an ALF space is the inverse face relation: *aNb* means that $a \in$SN(*b*) while *b* is a face of *a*.

There is an important particular case of ALF spaces, which is especially well appropriate for applications in computer imagery. It is characterized by a half-order relation between the elements of the space and by an additional feature: the dimension function *dim*(*a*), which assigns the smallest non-negative integer to each space element so that if $b \in$SON(*a*, *S*), then *dim*(*a*)≤*dim*(*b*). This kind of an LF space is called *abstract cell complex* or AC complex [14]. Its elements are called *cells*. If *dim*(*a*)=*k*, then *a* is called a *k-dimensional cell* or a *k-cell*. The dimension of a complex is the greatest dimension of its cells.

We have already introduced at the beginning of this Section the *bounding relation* $a<b$, which means $a \le b$ and $a \ne b$. It is irreflexive, asymmetric and transitive.

Dimensions of cells represent the half-order corresponding to the bounding relation. Let us call the sequence $a<b< \ldots <k$ of cells of a complex *C*, in which each cell bounds the next one, a *bounding path* from *a* to *k* in *C*. The number of cells in the sequence *minus one* is called the *length* of the bounding path.

**Definition DC:** (dimension of a cell) The dimension *dim*(*c*) of a cell *c* of a complex *C* is the length of the *longest bounding path* from any element of *C* to *c*.

This definition is in correspondence with the well-known notion of the dimension or height of an element of a partially ordered set [2a]

An example is shown in Fig. 4.1.

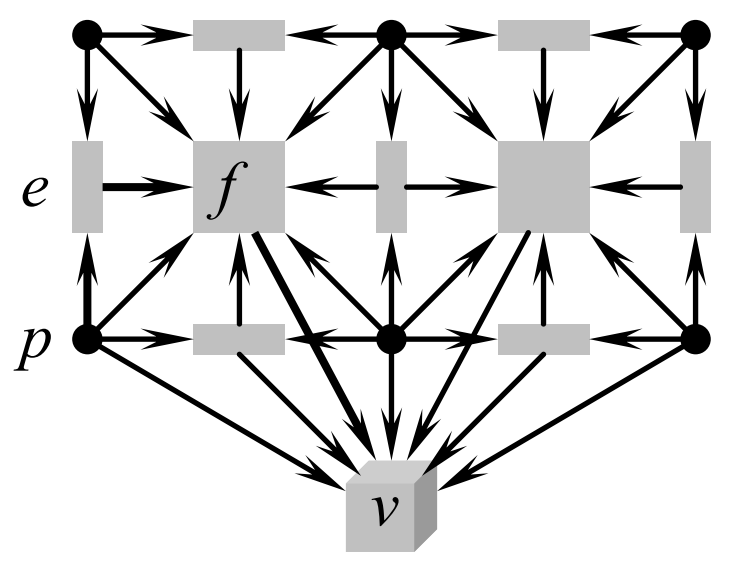

**Fig. 4.1**. A complex with bounding relations represented by arrows. An arrow points from *a* to *b* if *a* bounds *b*

Fig. 4.1 shows a complex and the bounding relations of its cells. The cell *p* is a minimum element of the space, its dimension is 0. One of the longest bounding paths from *p* to *v* is (*p*, *e*, *f*, *v*). Its length is 3, therefore *dim*(*v*)=3.

The dimension of the space elements is an important property. Using dimensions prevents one from some errors, which may occur when using an LF space without dimensions.

Thus, for example, it is possible to define a topology on the set $Z^n$ by defining the points with an odd value of $x_1+...+x_n$ as open and those with an even value of $x_1+...+x_n$ as closed [7, p. 199].

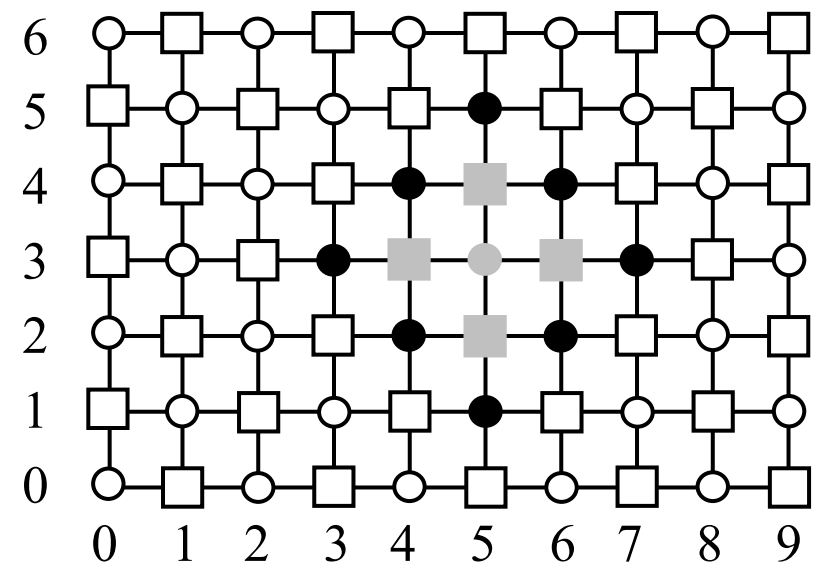

**Fig. 4.2.** The frontier (black disks) of a subset (gray squares and disks) under the one-dimensional topology assigned to $Z^2$

It is a topological space in the classical sense. However, it is a *one-dimensional* space and thus not appropriate for the *n*-dimensional set $Z^n$. Really, a closed point bounds 2*n* open points, which in turn bound no other points. Thus, the length of all bounding paths is equal to 1. It is easily seen that the frontier of a subset of this space is a *disconnected set* of closed points because the frontier contains no open points, while no two closed points are incident to each other. However, one expects from an *n*-dimensional space with *n*>1 that the frontier of its subset is connected. Fig. 4.2 shows the case of *n*=2. The squares represent the open points, circles and disks represent the closed ones. The neighborhood of a closed point contains adjacent squares, while the neighborhood of a square is the square itself. The neighborhoods of the closed points represented by black disks intersect both the foreground (gray) and the background (white). These are the only elements of the frontier and their set is disconnected. This error was made even by some experts because they have ignored the dimensions of space elements. The usage of dimensions of space elements is one of the advantages of the abstract cell complexes as compared with LF spaces without dimensions.

We shall use in the following sections a particular case of AC complexes known as Cartesian complexes [13, 15, 16]. Consider a space *S* whose elements compose a sequence:

$$S=(e_0, e_1, e_2, ..., e_{2m}) \tag{4.1}$$

Let each $e_i$ with an odd index *i* to bound the elements $e_{i-1}$ and $e_{i+1}$. Thus *S* becomes a one-dimensional ALF space. If we assign the dimension 0 to each closed and the dimension 1 to each open element, then *S* becomes a one-dimensional AC complex. Each closed cell is a face of two open cells. The indices in (4.1) are called *combinatorial coordinates* of cells in a one-dimensional space. (In former publications of the author [15, 16, 18] they were called "topological coordinates").

An AC complex of greater dimension is defined as the Cartesian product (also called the product set, set direct product, or cross product) of such one-dimensional complexes. A product complex is called a *Cartesian* AC complex [13]. The set of cells of an *n*-dimensional Cartesian AC complex $S^n$ is the Cartesian product of *n* sets of cells of one-dimensional AC complexes. These one-dimensional complexes are the *coordinate axes* of the *n*-dimensional space. They will be denoted by $A_i$, $i=1, 2, ..., n$. A cell *c* of the *n*-dimensional Cartesian AC complex $S^n$ is an *n-tuple* $c=(a_1, a_2, ..., a_n)$ of cells $a_i$ of the corresponding axes: $a_i \in A_i$. The cells $a_i$ are called *components* of *c*.

**Definition FRL:** (face relation) The face relation of the *n*-dimensional Cartesian complex $S^n$ is defined as follows: the *n-tuple* $(a_1, a_2, ..., a_n)$ is a face of another *n-tuple* $(b_1, b_2, ..., b_n)$ iff for all $i=1, 2,...n$ the cell $a_i$ is a face of $b_i$ in $A_i$.

Coordinates of a product cell *c* are defined by the vector whose coordinates are the coordinates of the components of *c* in their one-dimensional spaces.

**Theorem KC**: (*k*-dimensional cell) The dimension of a cell $c=(a_1, a_2,..., a_n)$ of an *n*-dimensional Cartesian complex *S* is equal to the number of its components $a_i$, $i=1, 2, ...n$; which are open in their axes.

It is easily seen, that the dimension of a product cell *c* is equal to the sum of dimensions of its components. Cells of the greatest dimension *n* are called the *principal cells*.

In the present paper we are interested to compare Cartesian complexes with standard 2D and 3D grids used in computer imagery. It is therefore appropriate to use here coordinates, which are *the same* for the principal cells of a Cartesian complex and for the pixels or voxels of a grid. These are the so called *semi-combinatorial coordinates*, which are integers for open components and half-integers for closed ones. If a cell has combinatorial coordinates $(x_1, x_2, ..., x_n)$, then its semi-combinatorial coordinates are $y_i=x_i/2$, $i=1, 2, ..., n$.

Cartesian complexes [13] are similar to spaces, which were independently published as Khalimsky spaces [6]. The difference consists only in the absence of the dimensions of space elements in Khalimsky spaces.

## 5 Previous Work

Since images are finite sets, it is natural to apply to them the axiomatic topology of locally finite spaces, especially that of AC complexes. There is, however, a difficulty, which has prevented until now the wide propagation of this theory among the researchers working in computer imagery. The difficulty arises due to the necessity to have space elements of different dimensions possessing different neighborhoods. Therefore attempts have been made to develop some kind of ersatz topology, which can be applied directly to sets of pixels or voxels. This is the concept of (*a*, *b*)-adjacencies [20, 10] often called the *graph-based approach to digital topology*. The ersatz topology is in concordance neither with the classical topology nor with that of LF spaces. However, it is widely spread in computer imagery since it is easy to understand and easy to apply.

As we have mentioned in the Introduction, many researcher in the field, including A. Rosenfeld, the author of the (*m*, *n*)-adjacencies, agree, that an adjacency relation is not consistent (or not “well-behaved”) if there is no topological space whose connectedness relation is analog to that of an image with that adjacency relation. Thus, to solve the problem posed by Yung Kong [9] and mentioned at the beginning of the Introduction, it is necessary to compare the connectedness relation in a digital image and in a corresponding topological space.

It is possible to try to define a topological space with the property that the connectedness of subsets is the same as in images provided with a particular (*m*, *n*)-adjacency. If the trial is successful, then this particular adjacency relation is “well-behaved”.

Both the graph-based and the axiomatic approach to digital topology have been discussed in the literature. In the textbook [7] both approaches are discussed “parallel”. This means that each problem is described two times, from both points of view. However, there is no such comparison of the connectedness as mentioned above. Klette has compared in his lecture notes [8] the connectedness of a (4, 8)-connected binary 2D image with that of a 2D complex and the (6, 26)-connected binary 3D image with that of a 3D complex. He has stated that when assigning the 0-dimensional and 1-dimensional cells to the foreground according to the “maximum-value rule” [14], then the connectedness in the image and in the complex are identical. However, he has not considered the general case.

Eckhardt and Latecki consider in their recent publication [4] three approaches to digital topology: the graph-theoretic approach, that of imbedding and the axiomatic approach. They consider as the axiomatic approach the Alexandroff topology [1]. They have demonstrated that the connectedness relation in a topological space, which is identical to that of a 2D image with the 8-adjacency, corresponds to a non-planar graph. But they have not considered pairs of adjacencies in 3D.

Kong and Rosenfeld [11] discuss the graph-based and the topological approach to digital topology. They consider a Euclidean complex *E* (without calling it by name) as a topological space, while the centers of the principal cells are points of $Z^n$. Thus $Z^n$ is identified with the set of the principal cells rather than with the whole set of cells of *E*, which is correct. The authors introduce the notion of “face-convex topological picture”. This is a complex with the following properties: given a subset *B* of the foreground (“black”) principal cells of *E*, the set *B** of all cells of *E* belonging to the foreground must satisfy

$$\mathrm{Int}(\mathrm{Cl}(B, S), S) \subseteq B^* \subseteq \mathrm{Cl}(B, S);$$

where Cl(*B*, *S*) is the closure and Int(*B*, *S*) the interior of *B* in *S*. (The definitions of these well-known notions are given below)

The authors have demonstrated that for any face-convex *n*-dimensional complex *EF*, *n*=2 or 3, there exists a digital picture with the specific adjacency relation depending on *B**, while *EF* is its “continuous analog”. The pair, in which one of the adjacencies is the "minimal" one, corresponding to $2\cdot n$ neighbors, and the other is the “maximal” one, corresponding to $3^n-1$ neighbors is a special case corresponding to *B**=Cl(*B*, *S*). However, they have not discussed the existence of analog complexes for digital images provided with all possible pairs of adjacencies.

In [7] the new idea of the switch-adjacency (s-adjacency) is suggested for 2D pictures. According to the s-adjacency exactly one of the two diagonally adjacent pairs of pixels is adjacent. It is possible to define one fixed pair as adjacent in all 2×2 blocks (Fig. 2.7 middle in [7]), or the choice of the pair may depend on the position of the 2×2 block of pixels. For example, one of the pairs is defined to be adjacent if the *y*-coordinate of the lower left-hand corner of the block is odd (Fig. 2.7 left in [7]). This two kinds of adjacency obviously correspond to a 6-adjacency applied to the rectangular grid as suggested in [12]. Another possibility is to "... let the s-adjacencies depend on the pixel values". This kind of s-adjacencies may be applied to grayscale pictures. This corresponds to the "membership rules" suggested in [14]. The authors have suggested no s-adjacency for the 3D case neither they have discussed the question of which adjacency pairs in 3D are "well-behaved".

## 6 Consistency of the (*m*, *n*)-Adjacencies from the Point of View of the Axiomatic Theory

We shall consider in this Section the (*m*, *n*)-adjacencies from the point of view of the axiomatic topology, find its limitations and suggest in Section 7 ways for topological investigations and for developing new topological algorithms, which are *directly applicable* to arrays of pixels or voxels.

Let us see, which adjacency pairs are consistent from the point of view of topology.

We shall need a universal notation for different adjacency relations in $Z^n$. The common notation by means of the number of adjacent points is inconvenient because these numbers depend on the dimension *n* of the space in a rather complicated way. We shall denote in what follows an adjacency relation between two points $P_1$ and $P_2$ of $Z^n$ by means of the *squared* Euclidean distance $d^2(P_1, P_2)=\sum(x_{1i}-x_{2i})^2$ between $P_1$ and $P_2$.

Let us call the set $\{x_1, x_2, ..., x_n\}$ of the coordinates of an element *e* of $Z^n$ or of a Cartesian complex the *coordinate set* of *e*.

**Definition CS:** (close) Two coordinate sets are called *close* to each other if the absolute values of all differences of their corresponding coordinates are less or equal to 1.

**Definition AD:** Two points $P_1$ and $P_2$ of $Z^n$ are called *a-adjacent* iff they are close to each other and $d^2(P_1, P_2)\le a$.

The value of *a* will be called the *index* of the corresponding adjacency relation. According to this notation the 4-adjacency in 2D becomes the 1-adjacency, the 8-adjacency becomes the 2-adjacency. In 3D the well-known 6-, 18- and 26-adjacencies become 1-, 2- and 3-adjacencies correspondingly. This notation may be easily used in $Z^n$ of any *n*.

Consider the quadruple ($Z^n$, *a*, *b*, *TR*), where *a* and *b* are indices of adjacency relations, *TR* is an *a*-connected subset of $Z^n$ called "foreground", while the subset $KR=Z^n-TR$ is the *b*-connected "background". We shall call, similarly as in [11], the quadruple ($Z^n$, *a*, *b*, *TR*) an *n-dimensional graph-based digital image*. On the other hand, we call the pair ($S^n$, *T*), where $S^n$ is an *n*-dimensional Cartesian complex partitioned into two subsets *T* and $S^n-T$, an *n-dimensional topological digital image*. Let $G^n$ be the set of the principal cells of $S^n$ and *M* be a one-to-one map between $Z^n$ and $G^n$ mapping each point *P* of $Z^n$ to the principal cell of $S^n$ having the same coordinates as *P*.

**Definition CL:** (closure) Let *t* and *T* be subsets of the space *S* such that $t\subseteq T\subseteq S$. Then the set containing with each cell $a\in t$ also all cells of *T*, which bound *a*, is called the *closure of t in T* and is denoted by Cl(*t*, *T*). The set *t*−Fr(*t*, *T*) is called the *interior of t in T* and is denoted by Int(*t*, *T*).

**Definition TA:** The *n*-dimensional topological image ($S^n$, *T*) is called the *topological analog* of the *n*-dimensional graph-based digital image ($Z^n$, *a*, *b*, *TR*) if the following conditions hold:

1. The closure Cl($M$($tr$), $T$) of $M$($tr$) of *any a-connected (a-disconnected) subset* $tr \subseteq TR$ is connected (correspondingly, disconnected).
2. The closure Cl($M$($cr$), $S^n - T$) of $M$($cr$) of *any* $b$-connected ($b$-disconnected) subset $cr \subseteq KR$ is connected (correspondingly, disconnected).

**Main Theorem:** There exists a topological analog ($S^n$, $T$) of the $n$-dimensional digital image $DI$=($Z^n$, $a$, $b$, $TR$) with an *arbitrary* subset $TR \subset Z^n$ iff $a \neq b$ and $a=n$ or $b=n$. A face-convex analog of $DI$ exits iff $a$=1, $b$=$n$ or $a$=$n$, $b$=1.

To prove the Theorem we need some definitions and lemmas.

**Definition IC:** (intermediate complex) If two principal cells $V_1$ and $V_2$ are close to each other, then the cell $C$ with the coordinates ($V_1$ + $V_2$)/2 is called the *intermediate cell* of $V_1$ and $V_2$. The closure $IC$=Cl({$C$}, $S^n$) is called the *intermediate complex* of $V_1$ and $V_2$. Note that $C$ is the greatest dimensional cell of $IC$.

**Lemma SC:** (small cell) If two principal cells $V$ and $W$ are close to each other, then each cell of the intermediate complex $IC$ is incident to both $V$ and $W$.

The proof is to be found in Appendix 1.

**Definition CS:** (corresponding subsets) A subset $T$ of $S^n$ is called *corresponding* to the subset $TR$ of $Z^n$ if there is a one-to-one correspondence between $Z^n$ and $S^n$, which maps the points of $TR$ onto the principal cells of $T$ and vice versa.

**Lemma NP:** (number of principal cells) A $k$-dimensional cell in an $n$-dimensional Cartesian complex is a face of $2^{(n-k)}$ principal cells.

We present the following proof here because it contains notions that we will need below.

**Proof of the Main Theorem:** Let $P_1$ and $P_2$ be two points of $TR$, which are $a$-adjacent and $V_1$ and $V_2$ be two principal cells of $T$ having the same coordinates as $P_1$ and $P_2$. The $a$-adjacency of $P_1$ and $P_2$ means according to Definition AD that $d^2(P_1, P_2) \leq a$ and hence $d^2(V_1, V_2) \leq a$. According to Lemma SC each cell $c$ of the intermediate complex $IC$=Cl({$C$}, $S^n$) of $V_1$ and $V_2$, where $C$=($V_1$+$V_2$)/2 is their intermediate cell, is incident to both $V_1$ and $V_2$. To make the set Cl({$V_1$}∪{$V_2$}, $T$) connected, it is necessary and sufficient to include any one cell of $IC$ in $T$. If, however, $d^2(V_1, V_2) > a$ then to make the set Cl({$V_1$}∪{$V_2$}, $T$) disconnected it is necessary to include the *whole intermediate complex IC* in the complement $K$=$S^n-T$. Similar conditions are guilty for points of $KR$ and the adjacency $b$ if we interchange $T$ and $K$. The demands of including a cell of $IC$ or the whole $IC$ in $T$ or in $K$ may cause contradictions.

Really, if the dimension $dim$($C$) of the intermediate cell $C$ is less than $n-1$ then, according to Lemma NP, there are at least four principal cells incident to $C$ and at least two pairs of principal cells, say ($V_1$, $V_2$) and ($W_1$, $W_2$), for which $C$ is its intermediate cell and $IC$=Cl({$C$}, $S^n$) is its intermediate complex. The distance $d^2(W_1, W_2)$ is equal to $d^2(V_1, V_2)$ since both of them are equal to $n-dim(C)$. If the $M$-preimages of $V_1$ and $V_2$ belong to $TR$ while the $M$-preimages of $W_1$ and $W_2$ belong to $KR$ and $d^2(W_1, W_2) \leq b$ then it is necessary and sufficient to include any cell of $IC$ in $S^n-T$. This does not contradict the demand that at least one cell of $IC$ be included in $T$ only if $IC$ contains more than one cell.

This is the case if $dim(C)>0$. Conclusions of these considerations are summed up in the following Table 1:

Table 1
Consistency of adjacency pairs

| | $d^2 \le b \Rightarrow c \in K$ | $d^2 > b \Rightarrow IC \subset T$ |
|---|---|---|
| $d^2 \le a \Rightarrow c \in T$ | OK iff $dim(C)>0$ | OK |
| $d^2 > a \Rightarrow IC \subset K$ | OK | contradiction |

"OK" means that there are no contradictions when including *c* or *IC* to *T* or *K*.

There are the following three cases of applying Table 1 to pairs of principal cells with different distances:

Case 1: $d^2(V_1, V_2)=n$; $dim(C)=n-d^2(V_1, V_2)=0$. Since $a,b \le n$; the expression $d^2 \le a$ in Table 1 can be replaced by $d^2=a$ and $d^2 \le b$ by $d^2=b$. Thus in this case there is no contradiction if *exactly one* of the indices *a* and *b* is equal to *n*; however, $a=b$ leads to a contradiction.

Case 2: $2 \le d^2(V_1, V_2) \le n-1$; $1 \le dim(C) \le n-2$. Table 1 can be used without changes. There is no contradiction if $1 \le dim(C) \le n-2$ and $d^2 \le a$ or $d^2 \le b$. Note that if $a=n$ or $b=n$, then the condition is always satisfied. Otherwise there may be a contradiction.

Case 3: $d^2(V_1, V_2)=1$; $dim(C)=n-1$. In this case the intermediate cell *C* is incident to a single pair of principal cells. Therefore, only one of the adjacencies *a* and *b* will be used and no contradiction can occur.

These cases can be summed up in the following table:

Table 2
Consistency depending on the dimension of the intermediate cell

| | | $b<n$ | | $b=n$ | |
|---|---|---|---|---|---|
| | | $b \ne a$ | $b=a$ | $b \ne a$ | $b=a$ |
| $a<n$ | $dim(C)=0$; case 1 | bad | bad | OK | – |
| | $0<dim(C)<n-1$; case 2 | variable | variable | OK | – |
| | $dim(C)=n-1$; case 3 | OK | OK | OK | – |
| | **all dimensions of *C*** | **bad** | **bad** | **OK** | – |
| $a=n$ | $dim(C)=0$; case 1 | OK | – | – | bad |
| | $0<dim(C)<n-1$; case 2 | OK | – | – | OK |
| | $dim(C)=n-1$; case 3 | OK | – | – | OK |
| | **all dimensions of *C*** | **OK** | – | – | **bad** |

Thus, the pair $(a, b)$ is consistent iff $a=n$ and $b<n$ or $b=n$ and $a<n$.

To obtain a *face-convex* topological analog it is necessary to include always the whole intermediate complex in one of the subsets $T$ or $K$. Table 1 must then be replaced by the following

Table 3
Face Convex consistency

| | $d^2 \le b \Rightarrow IC \subset K$ | $d^2 > b \Rightarrow IC \subset T$ |
|---|---|---|
| $d^2 \le a \Rightarrow IC \subset T$ | contradiction | OK |
| $d^2 > a \Rightarrow IC \subset K$ | OK | contradiction |

The above cases 1 and 3 can be applied to Table 3 without changes, but the case 2 must be changed:

Case 2: $2 \le d^2(V_1, V_2) \le n-1$; $1 \le dim(C) \le n-2$. According to Table 3, there is no contradiction if $a \ge d^2$ and $b < d^2$ or $b \ge d^2$ and $a < d^2$.

The most rigorous restriction for $a$ and $b$ takes place when $dim(C)=n-2$. Then $d^2=2$ and the consistency condition becomes

$$a \ge 2 \text{ and } b < 2 \text{ or } b \ge 2 \text{ and } a < 2. \qquad (6.1)$$

The condition of case 1 demands that exactly one of the indices $a$ and $b$ is equal to $n$. When combining this with (6.1) we arrive at the conclusion that a face-convex topological analog exists iff $a=n$ and $b=1$ or $b=n$ and $a=1$.

Let us show, that in the cases free of contradictions the desired connectedness of $T$ and $K$ is reached. A connected subset $tr \subseteq TR$ contains for any two points $P_1$ and $P_2$ of $tr$ an $a$-path lying completely in $tr$. The path is a sequence of pairwise $a$-adjacent points. Each such pair corresponds to a pair of principal cells of $S^n$ having the same coordinates. This pair is connected in Cl($t$, $T$) with $t=M(tr)$ due to the inclusion of cells lying between them to $T$ since Cl($t$, $T$) contains all cells of $T$ incident to the principal cells of $t$. Thus there is an incidence path in Cl($t$, $T$) between the principal cells corresponding to $P_1$ and $P_2$ and therefore $t$ is connected.

On the other hand, if a pair of points of $TR$ is $a$-disconnected, then the intermediate complex of the corresponding principal cells $V_1$ and $V_2$ of $T$ is included in $K$. In this case the set Cl($\{V_1\}$, $\{V_2\}$, $T$) is disconnected.

In a similar way one can demonstrate that Cl($k$, $S-T$) with $k=M(kr)$ is connected (disconnected) for any $b$-connected ($b$-disconnected) subset $kr \subseteq KR$. ■

Consider the examples of Fig. 6.1.

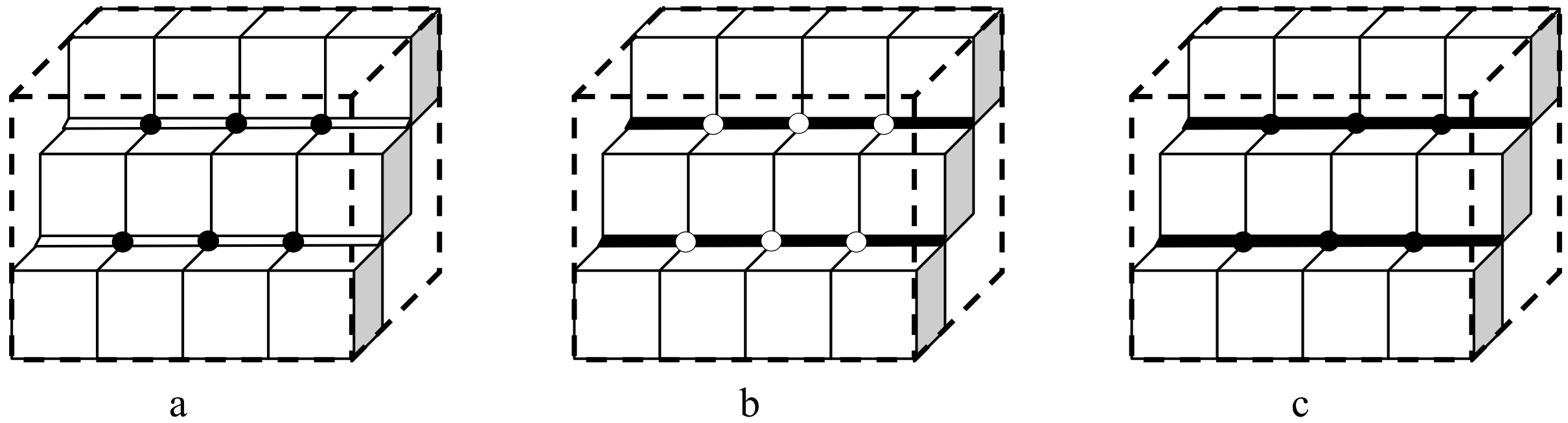

Fig. 6.1. Topological analogs of three graph-based digital images having the same set $TR$ and different adjacency pairs: (26,18) in (a), (18,26) in (b), and (6,26) in (c)

All three graph-based digital images are defined on the same subset of $Z^3$ containing 3×4×3 points. They have also the same subsets *TR* represented by the centers of visible cubes. Points of *KR* are the centers of the non-visible cubes. The 3-dimensional cells of $S^3$ are represented by cubes, the 2-cells by faces, the 1-cells by edges and the 0-cells by vertices of the cubes. The cells of *T* are the visible cubes, the faces lying between the visible cubes, the black edges and the black vertices. The white edges and vertices belong to $S^3-T$. The digital image of Fig. 6.1a has the (26,18)-adjacency. Each of the black vertices lies between two 26-adjacent but not 18-adjacent visible cubes. Thus the set *TR* corresponding to the set of visible cubes is 26-connected. Each of the white edges lies between two 18-adjacent non-visible cubes. Thus the set *KR* corresponding to the set of non-visible cubes is also connected. Other as in the case of the (8,8)-adjacency in 2D, this is no topological contradiction due to the fact that there are *three cells* in the set $IC^1$=Cl($\{p\}$, $S^3$)∩Cl($\{q\}$, $S^3$), where *p* and *q* are two 18-adjacent cubes. Thus some of the cells of $IC^1$ can by assigned to *T* and some other to $S^3-T$. This is impossible in the case of the (8,8)-adjacency because the set $IC^0$=Cl($\{p\}$, $S^2$)∩Cl($\{q\}$, $S^2$), where *p* and *q* are two 8-adjacent pixels, consist of a single cell.

The question may arise, whether a digital version of the Jordan theorem is guilty in such a space. The answer is “yes” for a topological space since a Jordan surface in a 3D ALF space is a 2D complex being a 2-dimensional combinatorial manifold [21]. It contains no 3-cells. The answer is “not” for a graph-based digital image since it is impossible to represent a 2D manifold in a 3D space as a graph-based digital image.

Fig. 6.1b shows the topological analog for the (18, 26)-adjacency. Also in this case both *TR* and *KR* are connected.

Fig. 6.1c shows the *face-convex* topological analog for the (26, 6)-adjacency. In this case *TR* is connected, but *KR* is disconnected.

According to the Main Theorem the adjacency pair (6,18) is not consistent. Some authors (e.g. [9]) assert that it is “well-behaved”. Let us demonstrate a counterexample. Consider a 3D image, in which the background is 18-connected and the 6-connected foreground consist of two cubes of *m*×*m*×*m* voxels, *m*≥3, which cubes have 2×2×2=8 common voxels at one corner (Fig. 6.2). Remove from each cube a smaller interior cube of (*m*–2)×(*m*–2)×(*m*–2) voxels. The remaining set *S* is a 6-connected simple “surface” according to the following often used definition (we write “surface” in quotation marks since a topological surface is a frontier of a solid 3D subset [17] and contains no 3-cells).

Let $N_b(p)$ be the *b*-adjacency of the point *p*.

**Definition SS:** (simple surface) The subset *B* of a graph-based digital picture *P*=($Z^n$, *a*, *b*, *B*) is a simple (*a*, *b*)-surface if it is *a*-connected and for each point *p*∈*B* the set $N_b(p)-\{p\}$ is a simple digital *b*-curve, i.e. a closed sequence *C* of *b*-connected points, in which any point is adjacent to exactly two other points of *C*.

It is easily seen, that the set *S* satisfies this definition. A simple surface must satisfy the Jordan Theorem and subdivide its complement into exactly two components.

However, the "surface" $S$ of Fig. 6.2 subdivides its 18-connected complement into three components, as it can be seen in Fig. 6.2b, which is a topological contradiction.

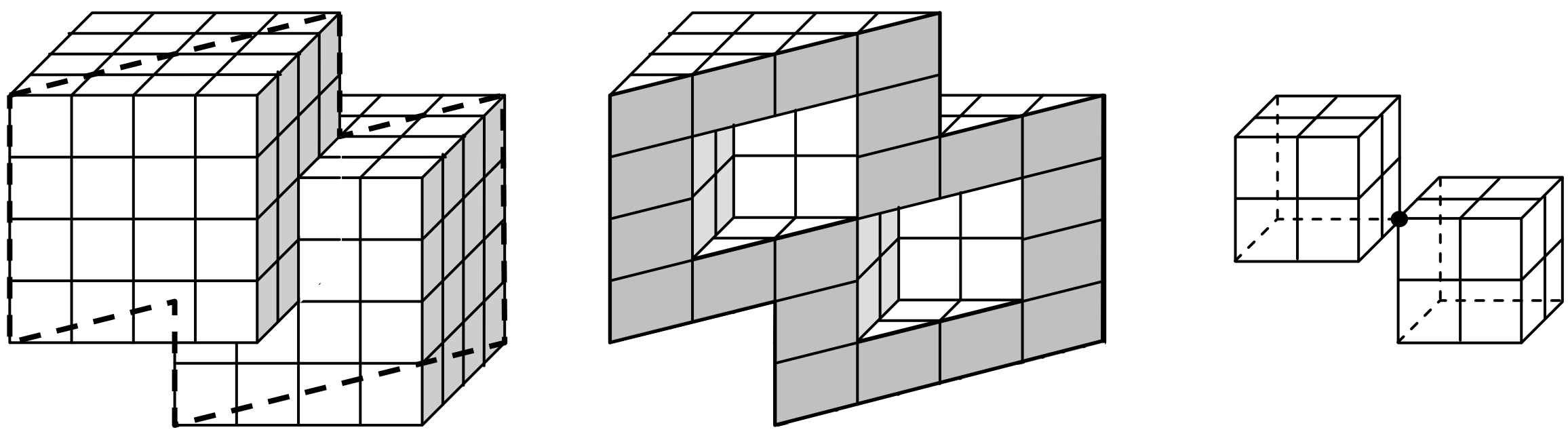

a b c

**Fig. 6.2.** Two hollow cubes (a), a cross section through the only common vertex of the interior cubes (b) and the interior cubes (c)

Yung Kong [9] considers also the pairs (12,12), (12,18) and (18,12) as "good pairs ... on the points of a 3-d face-centered cubic grid ... and (14,14) ... on the points of a 3-d body centered cubic grid". We have seen above that pairs ($a$, $b$) with $a$=$b$ are not good on $Z^3$. Why can pairs with $a$=$b$ be good on other grids?

There exist $n$-dimensional LF spaces with the property that any two principal cells that have a common incident cell of any dimension also have a common incident ($n$−1)-cell. We suggest to call them *completely connected spaces* or CC spaces. The notion of a strongly connected complex is well-known: an $n$-dimensional complex $S$ is called *strongly connected* if for any two principal cells $G_1$, $G_2$ of $S$ there exits a sequence of alternating mutually incident $n$- and ($n$−1)-cells of $S$, which sequence starts with $G_1$ and ends with $G_2$. *In an CC complex any connected subset is strongly connected.*

The simplest example is the hexagonal grid: any two hexagons that have a common vertex also have a common side. Also the 2D grid consisting of alternation squares and octagons is completely connected.

Completely connected spaces that correspond to tessellations of Euclidean 3-space by translates of a single polytope are of particular interest.

**Definition HCC:** An $n$-dimensional LF space, in which the principal cells are isomorphic to translates of a single polytope and any two principal cells, which have a common incident cell of any dimension also have a common incident ($n$−1)-cell, is called a *homogeneous completely connected space* or a HCC space.

As we have seen in the proof of the Main Theorem, if two principal cells have a common ($n$−1)-cell, then no contradiction can occur. In a HCC space *each close pair* has a common ($n$−1)-cell. Therefore in a HCC space the adjacency pair ($m$, $m$), where $m$ is the number of the principal cells incident to each principal cell, is consistent.

The reader can easily see that the hexagonal grid is the only two-dimensional HCC space. It is known [5] that the truncated octahedron (Fig. 7) is *the only* "primitive" parallelohedron, i. e., at each of its vertices exactly 4 (generally: $n$+1) tiles of a suitable face-to-face-tiling come together.

It has 14 faces. Parallelohedrons of this tiling are isomorphic to the principal cells of a 3D HCC space.

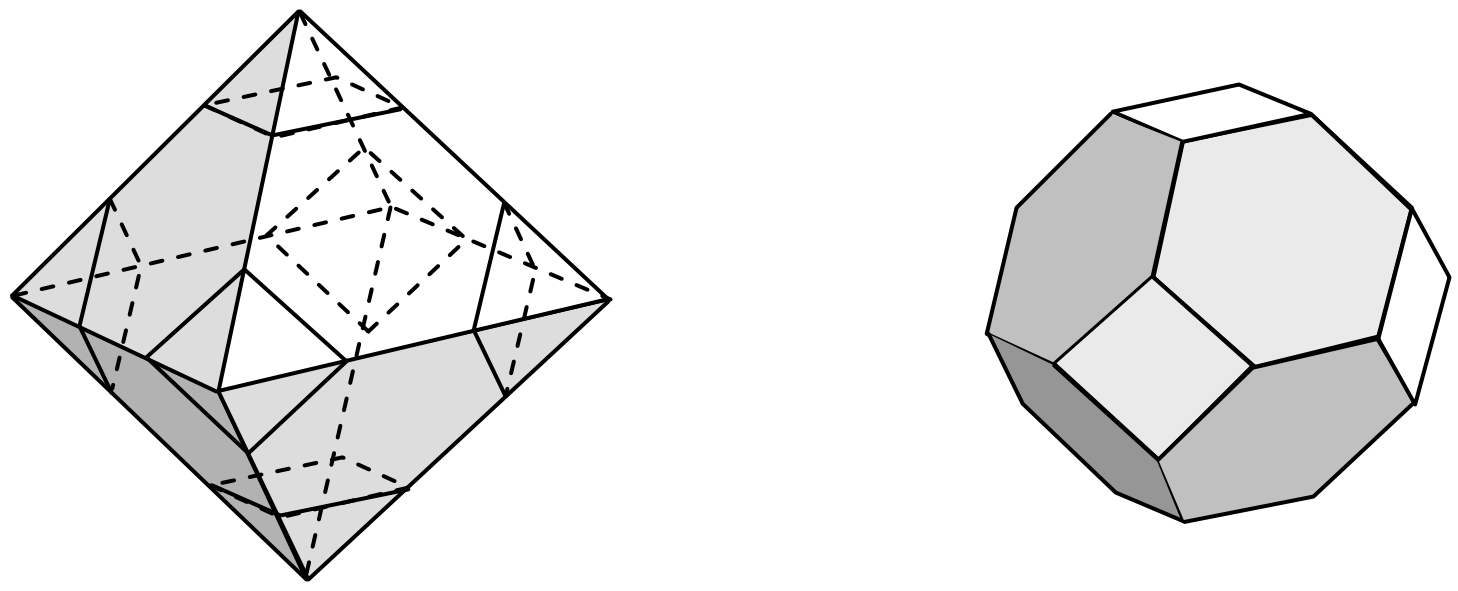

**Fig. 6.3**. The truncated octahedron

Since the truncated octahedron is the only primitive parallelohedron, other tesselations of the 3D space by translation do not compose HCC spaces. Therefore, in the tesselation by rhombic dodecahedrons, which corresponds to the above mentioned (12, 12) adjacency, two principal cells may have a 0-cell $C^0$ as the *single* common face. Thus a topological contradiction takes place when the first adjacency demands that $C^0$ belongs to the foreground while the second adjacency demands that it belongs to the background.

## 7 Recommendations for Applications

We have demonstrated that the majority of the adjacency pairs are topologically contradictory. Even the few consistent pairs, (4, 8) and (8, 4) in 2D and (6, 26) and (26, 6) in 3D, have important drawbacks:

a) They are only applicable in cases when there is only one subset of the space under consideration (and its complement). Thus they are not applicable, e.g., for colored images.

b) Even in the case of a black-and-white image the set created by an adjacency pair has a strange topological structure: thus, e.g. a 4-connected subset is full of holes due to the missing 0-cells. Consider e.g. an image with the (4, 8)-adjacency and the 4-connected dark subset of Fig. 7.1. It is topologically correct that the 0-cell between the pixels *d* and *f* belongs to the white subset making the pixels *d* and *f* disconnected and the pixels *c* and *g* connected.

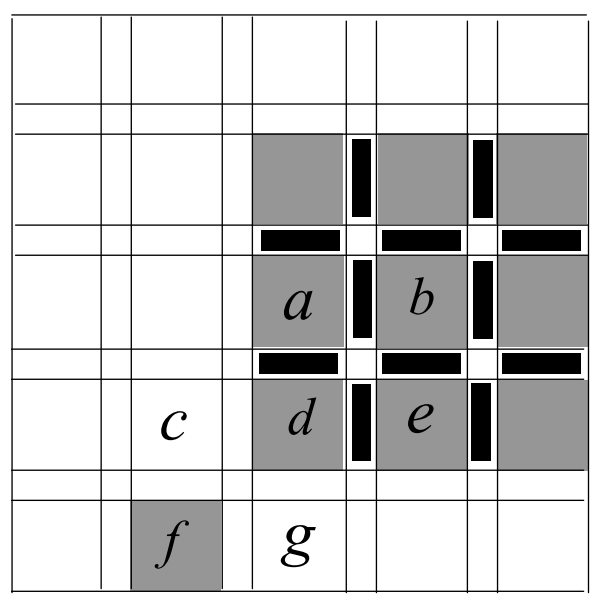

**Fig. 7.1.** Holes in a 4-connected subset

However, it is not necessary, that the 0-cell between the pixels *d* and *b* belongs to the white subset. This makes the greater component of the dark subset not simply connected. Each 0-cell incident to 4 pixels of *T* belongs to the frontier of the dark subset, and the frontier is disconnected, which is never the case for a simply connected subset of a topological space. This drawback does not affect the desired connectedness; it can, however, lead to topological contradictions when solving more complicated topological problems. For example, a 6-connected "surface" in 3D has a lot of tunnels due to missing 1-cells.

c) The frontiers of subsets under an adjacency relation are either not thin (compare Section 2) or the frontier of a subset is different from the frontier of its complement, which is a nonsense both from the point of view of topology and of the common sense.

d) The connectedness structure produced by an adjacency pair is no topological space at all since this structure is a property of a concrete subset of the space rather than that of the space. The connectedness changes when the subset changes.

What are the conclusions and the *recommendations for algorithm design*? We recommend not to use adjacency relations and to consider all topological and geometrical problems from the point of view of locally finite topological spaces (ALF spaces) or, even better, of complexes. The latter special case of an ALF space has some advantages due to the presence of the dimensions of cells. The dimensions of cells make the work with the topological space easier and more illustrative.

The usual objections against the use of complexes are the following:

a) Why should we use cells of lower dimension, which we don't see on a display?

b) When using complexes, one needs much more memory space: 4 times or 8 times more in the 2D or in the 3D case correspondingly.

The objection because of visibility is not pertinent since the visibility has nothing to do with topology. For example, we all use in our work with 3D images the voxels. However, the *voxels are not visible* on the displays: what we see are the *faces of voxels*, i.e. the 2-cells, while voxels are 3-cells. Thus, e.g. the software "OpenGL", which is widely spread for representing 3D scenes, works only with faces of polyhedrons like triangles, squares, polygons and does not use three-dimensional bodies at all.

The second objection is pertinent only if one would try to allocate memory space for cells of all dimensions, which is almost never necessary. Cell complexes are a *means for thinking* about geometrical and topological problems rather than for data saving. It is possible to work with complexes, while saving in the memory only certain values assigned to the principal cells, like colors for the pixels, or densities for the voxels. Cells of lower dimensions are present as some kind of *virtual objects* only. Algorithms of this kind are described in [18, 19].

In cases when the membership of cells of lower dimensions is important it can be defined by a *"flat rate" face membership*, i.e., a rule specifying the set membership of each cell as a function of the membership of the incident principal cells.

Consider the example of Fig. 7.2. Suppose, it is necessary to define the connectedness in such a way that both the white and the black "V" are connected. This is obviously impossible under any adjacency relation, neither under the s-adjacency [7]. The aim can be achieved by means of the following rule [14] for assigning membership labels to 1- and 0-cells.

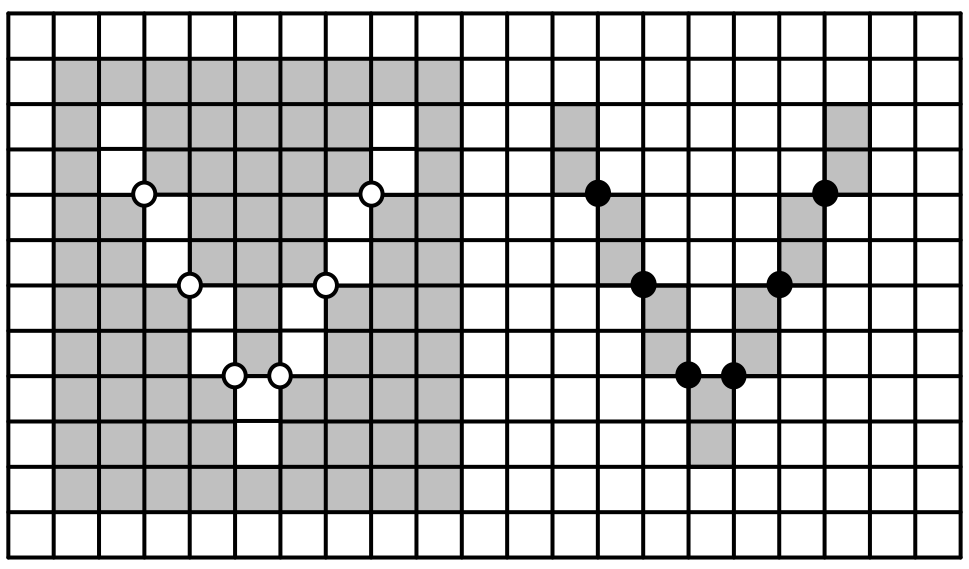

**Fig. 7.2.** A white and a black V-shaped regions in one image, both connected due to applying the EquNaLi rule

**EquNaLi Rule:**

A 1-cell gets the greatest label of the two incident pixels.

The label of a 0-cell $c^0$ is defined as follows:

If SON($c^0$) contains exactly one diagonal pair of pixels with equivalent ("Equ") labels, then $c^0$ gets this label. If there are two such pairs but only one of them belongs to a narrow ("Na") stripe, then $c^0$ gets the label of the narrow stripe. Otherwise $c^0$ gets the maximum, i.e. the lighter ("Li") label of the pixels of SON($c^0$).

The latter case corresponds to the cases when SON($c^0$) contains no diagonal pair with equivalent labels or it contains two such pairs and both of them belong to a narrow stripe.

To decide, whether a diagonal pair in SON($c^0$) belongs to a narrow stripe it is necessary to scan an array of 4×4 pixels with $c^0$ in the middle and to count the labels corresponding to both diagonals. The smaller count indicates the narrow stripe. Examples of other efficient membership rules may be found in [14].

Another important problem is that of using completely connected spaces. It is not correct to think that we need some special scanners or other special hardware to work with a hexagonal 2D space or with a 3D space tessellated by truncated octahedrons. It is possible to use as ever the standard orthogonal grids. The only thing, which must be changed is the definition of the incidence and hence that of the connectedness.

Consider the example of the hexagonal 2D space.

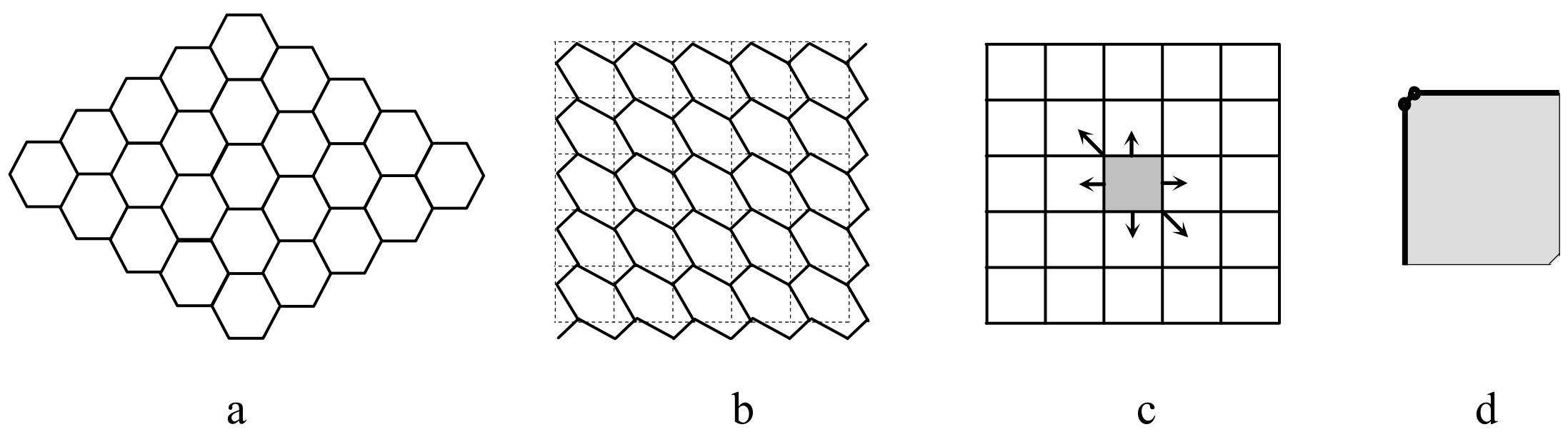

**Fig. 7.3.** A hexagonal grid (a), a transformed rectangular grid (b), a rectangular grid with 6-adjacency (c) and the "virtual" cells, transforming each pixel to a hexagon

Fig. 7.3a shows a hexagonal grid. In Fig. 7.3b we see the transformed rectangular grid, in which the upper left corner of each rectangular pixel is replaced by a slanted edge. The original rectangular grid is shown in Fig. 7.3c. The arrows show the only six pixels, which should be considered as adjacent to the shaded pixel. If it is necessary to take cells of lower dimensions in consideration, then two 0-cells and three 1-cells must be assigned to each pixel, as shown in Fig. 7.3d. Thus two *virtual* 0-cells and one *virtual* 1-cell correspond to each vertex of the grid. The remaining two 1-cells correspond to the edges of the grid. All five faces assigned to a pixel have the same coordinates as the pixel. They can be distinguished by their types only. If it is necessary to assign a label to each 0-cell and to each 1-cell then 5 bits of the memory word assigned to the pixel can be used for these cells. The rest of the memory word can be used for the label or for the color of the pixel. This data structure and some algorithms working in a *virtual* hexagonal grid are described in [12]. It should be stressed once more, that the 6-adjacency realized by this data structure may be used for multicolored images and causes no topological paradoxes.

A similar data structure can be easily developed for the 3D standard grid to realize the 14-adjacency corresponding to the only 3D HCC space.

## Conclusion

We have suggested a new set of Axioms of digital topology and have demonstrated, that they define a space, which is a particular case of a classical topological space. This fact serves as an explanation, why the key object of the classical topology, the system of open subsets, must have the properties formulated as the classical axioms: these properties are necessary (but not sufficient) to define the connectedness through neighborhoods and the frontier of a subset $T$ as a thin subset, which is the same for a subset $T$ and for its complement. Thus we can now answer the second question, raised in the Introduction: yes, it is possible to find such axioms, however, they lead to a concept of a space, which is a special case of a classical topological space.

The paper demonstrates how the ($a$, $b$)-adjacency relations commonly used in computer imagery can be brought into accordance with the connectedness of a topological space. It was demonstrated that in spaces of any dimension $n$ only those pairs ($a$, $b$) of adjacencies are consistent, in which exactly one of the adjacencies is the "maximal" one corresponding to $3^n-1$ neighbors. Even the consistent pairs have important limitations: they are not applicable for multicolored images and they cannot correctly represent topological properties of subsets. We suggest to use instead the concept of ALF spaces, while considering space elements of lower dimensions as "virtual" objects, which need not be saved in the memory. This attitude makes it possible to apply consistent topological definitions and algorithms to images represented in standard grids. The notion of homogeneous completely connected spaces is introduced, and it has been demonstrated, that there is only one such space of dimension 3; it is isomorphic to the tessellation by truncated octahedrons, each of which has 14 faces.

The paper shows, how locally finite spaces, especially cell complexes and completely connected spaces, can be applied to computer imagery, while using standard orthogonal grids.

## Acknowledgements

The author is pleased to acknowledge many useful discussions with Boris Flach, Siegfried Fuchs and Henrik Schulz on parts of the work presented here.

## References

[1] P. Alexandroff, *Diskrete Topologische Raeume*, Matematicheskii Sbornik 2, 1937

[2] P. Alexandroff and H. Hopf, *Topologie I*, Springer, 1935

[2a] G. Birkhoff, Lattice Theory, American Mathematical Society, 1961.

[3] U. Eckhardt and L. Latecki, Topologies for the digital spaces $Z^2$ and $Z^3$. *Computer Vision Image Understanding, 90*, 2003, pp. 295-312

[4] U. Eckhardt and L. Latecki, *Digital Topology*, http//cosmic.rrz.uni-hamburg.de/docs/wcat/mathematik/eckhardt/eck0001/eck0001.pdf

[5] E.S. Fedorov, "The Symmetry of Regular Systems of Figures", *Zap. Miner. Obshch.* 21 (1885), 1-279. Reprinted, Moscow, Izdat. Akad. Nauk SSSR, 1953. (in Russian.)

[6] E. Khalimsky, R. Kopperman and P.R. Meyer, "Computer Graphics and Connected Topologies on Finite Ordered Sets", *Topology and Applications*, 36, 1990, pp. 1-17

[7] R. Klette and A. Rosenfeld, *Digital Geometry*, Morgan Kaufmann, 2004.

[8] R. Klette,. *Lecture Notes*, www citr.auckland.ac.nz/~rklette/Books/MK2004/pdf-LectureNotes/

[9] T.Y. Kong, *Can 3-D Digital Topology be Based on Axiomatically Defined Digital Spaces?* http //citr.auckland.ac.nz/dgt, Open Problems, 2001

[10] T.Y. Kong and A. Rosenfeld, "Digital Topology. Introduction and Survey", *Computer Vision, Graphics and Image Processing*, vol. 48, 1989, pp. 357-393

[11] T.Y. Kong and A. Rosenfeld, "Digital Topology, a Comparison of the Graph-based and Topological Approaches", in G.M. Reed, A.W. Ronscoe, and R.F. Wachter (ed), *Topology and Category Theory in Computer Science,* Oxford University press, Oxford, U.K, 1991, pp. 273-289.

[12] V. Kovalevsky, "Discrete Topology and Contour Definition", *Pattern Recognition Letters*, vol. 2, 1984, pp. 281-288

[13] V. Kovalevsky, "On the Topology of Digital Spaces", *Proceedings of the seminar "Digital Image Processing"* in Fleeth/Mirrow, Germany. Technical University Dresden 1986, pp. 1-16

[14] V. Kovalevsky, "Finite Topology as Applied to Image Analysis", *Computer Vision, Graphics and Image Processing*, vol. 46, No. 2, 1989, pp. 141-161

[15] V. Kovalevsky, "Finite Topology and Image Analysis", in P. Hawkes (ed.), *Advances in Electronics and Electron Physics*, Academic Press, v. 84, 1992, pp. 197-259

[16] V. Kovalevsky, "Digital Geometry Based on the Topology of Abstract Cell Complexes", in *Proceedings of the Third International Colloquium "Discrete Geometry for Computer Imagery"*, University of Strasbourg, September 20-21, 1993, pp. 259-284

[17] V. Kovalevsky, "A New Concept for Digital Geometry", in *Shape in Picture*, Springer-Verlag 1994, pp. 37-51

[18] V. Kovalevsky, "Algorithms and Data Structures for Computer Topology", in G. Bertrand et all (eds.), LNCS 2243, Springer-Verlag 2001, pp. 37-58

[19] V. Kovalevsky, "Algorithms in Digital Geometry Based on Cellular Topology", in R. Klette. and J. Zunic. (eds), LNCS 3322, Springer Verlag 2004, pp. 366-393

[20] A. Rosenfeld and A. C. Kak, *Digital Picture Processing*, Academic Press, 1976

[21] J. Stillwell, *Classical Topology and Combinatorial Group Theory*, Springer, 1995

**Appendix 1: The Proofs**

**Proof of Theorem TF:** According to Definition TF the frontier Fr(*T*, *S*) is not thin if it contains an opponent pair (*a*, *b*). According to Definition OT (Section 2), the relations $a \neq b$, *aNb* and *bNa* must be true for this pair. If *N* is antisymmetric, then the relations cannot be true simultaneously. Therefore no opponent pair can exist.

On the other hand, if *N* is not antisymmetric, then there exist $a,b \in S$ such that $a \neq b$, *aNb* and *bNa*. Then *a* belongs to the frontier Fr({*a*}, *S*) since SN(*a*) contains an element of {*a*} and an element *b* which is not in {a}. In the same way we can see that *b* also belongs to Fr({*a*}, *S*). Therefore, *a* and *b* are opponents and the frontier Fr({*a*}, *S*) is not thin. ■

**Proof of Lemma SI:** Suppose *T* is open, $a \in T$ and there exists $b \in$ SN(*a*), which belongs to the complement of *T*. Then according to Definition FR *a* belongs to Fr(*T*, *S*), which contradicts the Definition OP. Thus, if a subset *T* is open, then it contains the SNs of all its elements. Now suppose that a subset *T* contains the SNs of all its elements and is not open. This means that it contains an element, say $b \in T$, which belongs to Fr(*T*, *S*). According to Definition FR SN(*b*) must intersect the complement $S-T$. Thus *T* does not contain SN(*b*). The contradiction proofs the Lemma.■

**Proof of Theorem OS:** The subset ∅ contains no elements at all and thus the conditions of Lemma SI are fulfilled in a trivial way. Thus according to Definition OP the subset ∅ is open in *S*. The frontier Fr(*S*, *S*) of the entire space is obviously empty since the complement of *S* is empty. Thus *S* contains no elements of Fr(*S*, *S*) and therefore *S* is open in itself.

Let the open subsets be $O_i$, $i$=1, 2, .... Consider an element *a* of the union $U = \cup O_i$, $i$=1, 2, ...; that belongs to $O_k$, where $k \geq 1$. Since $O_k$ is open, it contains according to Lemma SI all elements of *S* that are in SN(*a*). So does the set *U*. Similar arguments are valid for each value of *k*. Thus *U* contains all elements of *S* that are in SN(*a*) where *a* is any element of *U*. According to Lemma SI *U* is open in *S*. Now we shall deduce a *stronger* version of the third axiom C3, which is well-known to be valid for Alexandroff spaces [1]. Consider an element *a* of the intersection $I = \cap O_i$, $i$=1, 2, .... It belongs to each of $O_i$. Since each $O_i$ is open, it contains all elements of *S* that are in SN(*a*). So does the set *I* since it contains all elements of *S*, which are contained in each $O_i$. According to Lemma SI the set *I* is open in *S*. The third classical axiom C3 follows immediately.■

**Proof of Lemma MM:** As we know, the bounding relation *B* is irreflexive and asymmetric. If we suppose, that *B* is transitive, then it is a half-order. Thus we can write in this case $a<b$ instead of *aBb*. Remember that *aBb* means $a \neq b$ and $b \in$ SN(*a*). Consider three elements *a*, *b* and *c* such that *aBb* and *bBc*. Since *B* is transitive, the relation *aBc* follows, which means that $c \in$ SN(*a*) and SN(*b*)⊆SN(*a*) since *c* may be any element of SN(*b*). However, $a \notin$ SN(*b*); thus SN(*a*) contains at least one element more than SN(*b*).

Consider a sequence of elements, in which each element bounds the next one:

$$...<a<b<c<d<e\ ...$$

Since the SN of each next element of the sequence contains less elements than the previous one, the sequence must be finite. Otherwise the SN of the left most element of the sequence would contain an unbounded number of elements, which is impossible since the space *S* is locally finite. Thus there exists the left most element, which is bounded by no other elements of *S*. This is a minimum element of *S*.

The sequence also cannot be infinitely continued to the right hand side. Really, since the number of elements of SN(*a*) is finite and the SN of each next element contains less elements than SN($a$), there must be in the sequence the right most element. Its SN contains only one element. This element bounds no other elements. This is a maximum element of $S$. ■

**Proof of Lemma NM:** If $a \in Fr(T, S)$, then SN($a$) must intersect both $T$ and $S-T$. Thus SN($a$) must contain at least two elements, one belonging to $T$ and another belonging to $S-T$. Let $b \neq a$ and $b \in SN(a)$. This means $a<b$. Therefore $a$ is no maximum element.

Let us prove now the second assertion. Let $b \in SN(a)$, which means $a<b$, and suppose that $b$ is not a maximum element of $S$. Then there must be in $S$ an element $c$ such that $b<c$, which means $c \in SN(b)$. Since $S$ is an LF space, the sequence $a<b<c<...$ must finish at a maximum element. Since $B$ is transitive, all elements of the sequence are in SN($a$). Therefore the maximum element belongs to SN($a$). ■

**Proof of Theorem TR:** Let $T$ be a subset of $S$ and $F$=Fr($T$, $S$) be the frontier of $T$. We must prove, that

1. if the neighborhood relation is transitive, then Fr($F$, $S$)=$F$ is fulfilled for all subsets $T \subset S$ and
2. if the condition Fr($F$, $S$)=$F$ is fulfilled for all subsets $T \subset S$, then the neighborhood relation is transitive.

Let us first prove the assertion 1. Let $a \in F$. It follows from $a \in SN(a)$ that $SN(a) \cap F \neq \emptyset$. According to Lemma MM, the set SN($a$) contains a maximum element of $S$ which, according to Lemma NM, never belongs to a frontier. Thus $SN(a) \cap (S-F) \neq \emptyset$ and the conditions for $a \in Fr(F, S)$ are fulfilled and any element of $F$ belongs to Fr($F$, $S$).

Now let $b \in Fr(F, S)$ which means $SN(b) \cap F \neq \emptyset$ and $SN(b) \cap (S-F) \neq \emptyset$. The latter condition is always fulfilled since SN($b$) contains a maximum element. The first condition means that there is in SN($b$) an element $c \in F$. Transitivity of $N$ means that any element $d \in SN(c)$ belongs to SN($b$). Thus $SN(c) \subseteq SN(b)$. Since SN($c$) intersects both $T$ and $S-T$, so does SN($b$). Therefore, each element of Fr($F$, $S$) belongs to $F$ and the assertion 1 is true.

To prove assertion 2 let $N$ be non-transitive. Then there exist distinct elements $a$, $b$, $c \in S$ such that $b \in SN(a)$, $c \in SN(b)$, however, $c \notin SN(a)$. Consider the singleton $\{c\}$. It follows from $c \notin SN(a)$ that $a \notin Fr(\{c\}, S)$ since $SN(a) \cap \{c\} = \emptyset$. On the other hand, $a \in Fr(Fr(\{c\}, S), S)$. Really, $b$ bounds $c$ and thus $b \in Fr(\{c\}, S)$. It follows from $b \in SN(a)$ that $SN(a) \cap Fr(\{c\}, S) \neq \emptyset$. Also $SN(a) \cap (S-Fr(\{c\}, S))$ is not empty since $a \notin Fr(\{c\}, S)$. Thus, both assertion 1 and 2 are proved. ■

**Proof of Corollary NO:** Let $a$ be an element of $S$ and SN($a$) its smallest neighborhood. Suppose SN($a$) is not open in $S$ according to Definition OP (Section 3). Then SN($a$) must contain at least one element $b \in Fr(SN(a), S)$. This means that SN($b$) contains at least one element $c \notin SN(a)$. Since $a<b$, $b<c$ and the bounding relation is transitive, the relation $a<c$ is true, which means that $c \in SN(a)$. The contradiction proves that SN($a$) is open according to Definition OP. However, according to Theorem OS (Section 3), subsets open according to Definition OP are also open in the classical sense. To prove that SN($a$) is the smallest open subset of $S$ containing $a$ let us remove any element $b \neq a$ from SN(a). The reduced set does not contain the smallest neighborhood of $a$. Thus, according to Lemma SI (Section 3) it is not open. This proves the Corollary. ■

**Proof of Corollary T0**: Consider any two space elements $a$ and $b$. If they are not incident to each other, then $a \notin SN(b)$ and $b \notin SN(a)$. In this case each of the sets SN($a$) and SN($b$), being open according to Corollary NO, satisfies the condition of Axiom $T_0$. If, however, they are incident, then either $b \in SN(a)$ or $a \in SN(b)$. Since the neighborhood relation is antisymmetric, the condition $b \in SN(a)$ implies $a \notin SN(b)$. In this case the open subset SN($b$) satisfies the condition of Axiom $T_0$. In the case when $a \in SN(b)$ this is the set SN($a$).■

**Proof of Theorem KC:** Let $c=(a_1, a_2,..., a_n)$ be a cell of $S$ having exactly $k$ open components. Then we can find a sequence of cells beginning with $c$, such that each next cell in the sequence has exactly one component, which bounds one of the open components of the previous cell. Thus, when changing from one cell of the sequence to the next one, exactly one open component will be replaced by a closed component. The last cell in the sequence is a 0-cell having all components closed. The sequence is a bounding path and its length is $k$. Therefore according to Definition DC, the dimension of $c$ is equal to $k$. ■

**Proof of Lemma SC:** According to Definitions FC and FRL (Section 4) a cell $c_1$ of a Cartesian complex is a face of the cell $c_2$ if each component of $c_1$ is a face of the corresponding component of $c_2$. In a Cartesian complex with semi-combinatorial coordinates (Section 4) a component with a half-integer coordinate $x$ is a face of two components with integer coordinates $x \pm 1/2$. A component with an integer coordinate is a (non-proper) face of itself only. Therefore, a cell $c$ is a face of a principal cell $V$ iff $|c_i-V_i| \leq 1/2$ for all $i \in (1, n)$ ($c$, $V$ and $W$ are considered here as $n$-dimensional vectors). A cell $c$ is a face of two close principal cells $V$ and $W$ if for all $i \in (1, n)$ either $V_i=W_i$ and $|c_i-V_i| \leq 1/2$ or $V_i \neq W_i$ and $c_i=(V_i+W_i)/2$. Thus there are $3^{(n-k)}$ common faces of two close principal cells $V$ and $W$, where $k=d^2(V, W)$ is the number of different coordinates of $V$ and $W$. A face $c$ of $V$ is incident to $V$. According to Theorem KC (Section 4) the dimension of $c$ is equal to the number of its integer coordinates. The cell with the greatest dimension among the common faces of $V$ and $W$ is obviously the intermediate cell $C=(V + W)/2$ since $C_i$ are integers *for all i*, for which $V_i=W_i$ and $C_i$ can not be integer if $V_i \neq W_i$. It is easily seen that all common faces of $V$ and $W$ are faces of $C$. Therefore, they compose the closure of $C$. ■

**Proof of Lemma NP:** According to Theorem KC (Section 4) a $k$-cell $c^k$ has $k$ integer coordinates. When adding to each of the remaining $n-k$ coordinates $\pm 1/2$ one gets a coordinate set containing $n$ integers, i.e. a coordinate set of a principal cell incident to $c^k$. There are $2^{(n-k)}$ possibilities of adding. ■

## Appendix 2: List of Abbreviations

2D; 3D; *n*D – two-dimensional; three-dimensional; *n*-dimensional;
AC complex – abstract cell complex; introduced in Section 4
ALF space – locally finite space satisfying our Axioms; introduced in Section 1
CC space – completely connected space; introduced at the end of Section 6
Cl(*B*, *S*) –closure of the subset *B* of the space *S*; Definition CL in Section 6
$d^2(P_1, P_2)$ – squared Euclidean distance between the points $P_1$ and $P_2$; introduced in Section 6
*dim*(*a*) – dimension function. Introduced in Section 4; Definition DC in Section 4
Fr(*T*, *S*) – frontier of the subset $T \subseteq S$ of the space *S*; Definition FR in Section 2
$G^n$ – the set of the principal cells; principal cells are the cells of the greatest dimension
HCC space – homogeneous completely connected space; Definition HCC in Section 6
iff – if and only if (a common notation in mathematical literature)
Int(*B*, *S*) – interior of the subset *B* of the space *S*; Definition CL in Section 6
LF space – locally finite space; Definition LFS in Section 2
*N* – neighborhood relation; Definition NR in Section 2
SN(*e*) – smallest neighborhood of the space element *e*; introduced in Section 2
SON(*a*, *S*) – smallest open neighborhood of space element *a* in the space *S*; Definition SON in Section 4
$Z^n$ – Cartesian product of *n* sets of all integers (a common notation in mathematical literature)
$S^n$ – *n*-dimensional Cartesian complex or *n*-dimensional space; introduced in Section 6